\documentclass[11pt]{article}
\usepackage{amsfonts, amsmath, amssymb, amsthm, float, fancyhdr}
\usepackage[pdftex]{graphicx}
\usepackage[a4paper,vmargin={20mm,20mm},hmargin={20mm,20mm}]{geometry}
\usepackage{enumerate}

\everymath{\displaystyle}
\DeclareGraphicsExtensions{.pdf,.png,.jpg}

\begin{document}
\begin{center}
{\huge{Lipschitz Functions on Sparse Graphs}}

\bigskip
Samuel Korsky, Tahsin Saffat, Dhroova Aiylam

\bigskip
October 2, 2023

\bigskip
{\bf{Abstract}}
\begin{changemargin}{2cm}{2cm} 
In this work we attempt to count the number of integer-valued $h$-Lipschitz functions (functions that change by at most $h$ along edges) on two classes of sparse graphs; grid graphs $L_{m,n}$, and sparse random graphs $G(n,d/n)$. We find that for all  $n$-vertex graphs $G$ with $k$ connected components, the number of such functions grows as $(ch)^{n - k}$ for some $1 \le c \le 2$. In particular, letting $\alpha \approx 1.16234$ be the largest solution to $\tan{(1/x)} = x$, we prove that as $n \to \infty$
\begin{itemize}
\item
$c = \alpha\sqrt{2} \approx 1.6438$ when $G = L_{2,n}$
\item
$1.351 \approx \alpha^2 \le c \le \arctan{(3/4)}^{-1} \approx 1.554$ when $G = L_{n,n}$
\item
$ 1 + \frac{1}{2d} + O\left(\frac{1}{d^2}\right) \le c \le 1 + \frac{4\ln^2{d}}{d} + O\left(\frac{1}{d}\right)$  (w.h.p.) when $G = G(n, d/n)$
\end{itemize}
\end{changemargin} 
\end{center}
\section{Introduction}
An $h$-Lipschitz function on a connected graph $G$ is an integer-valued function on the vertices of $G$ which changes by at most $h$ between adjacent vertices, and where the function's value at a fixed vertex in $G$ is $0$. In this work we seek to count the number of such functions in the setting where $h$ is large, for various classes of sparse graphs; in particular, for grid graphs and sparse Erd\H{o}s-Renyi random graphs. 

\bigskip
\noindent
Our motivation for this work comes from the extensive literature pertaining to $h$-Lipschitz functions [1, 2, 3] and the related concept of $\mathbb{Z}$-homomorphisms (integer functions on the vertices of $G$ which change by exactly one along edges) [4, 5]. In particular, the study of $h$-Lipschitz functions on $d$-regular expanders [1] is closely related to our work on random graphs, and previous work on $h$-Lipschitz functions on tori [2] and $\mathbb{Z}$-homomorphisms on lattices [4, 5] partially inspires our investigation on grid graphs. 

\bigskip
\noindent
The principal difference between our work and the existing literature is that we study $h$-Lipschitz functions for large $h$ (indeed, we take $h \to \infty$ for the majority of our results). In previous work [1, 2, 3], $h$ is generally bounded as a function of the average or maximum degree of the graph in question. And $\mathbb{Z}$-homomorphisms are more different still, with their closest equivalent being $h$-Lipschitz functions for $h = 1$. 

\bigskip
\noindent
Note that (as mentioned in [1]) our work also has potential applications to surface models in statistical mechanics; counting the number of $h$-Lipschitz functions on a large grid graph can be interpreted as finding the entropy of surface energy models (satisfying certain constraints) at high temperature, where the probability distribution over valid energy configurations is uniform. 

\subsection{Definitions}

We begin with a formal definition: 

\bigskip
\noindent
{\bf{Definition 1.1.1}}. For a given connected graph $G$ and a fixed vertex $v \in V(G)$, let $\text{Lip}_{v}(G; h)$ denote the set of $h$-Lipschitz functions from the vertices of $G$ to $\mathbb{Z}$ that send $v$ to $0$; that is, functions $f: V(G) \rightarrow \mathbb{Z}$ such that $f(v) = 0$ and $|f(x) - f(y)| \le h$ for every $\{x, y\} \in E(G)$. It is clear that $\left\lvert\text{Lip}_{v}(G; h)\right\rvert$ remains the same regardless of choice of $v$, so we will generally eschew the ``$v$'' subscript as we are primarily interested in the cardinality of the set. 

\bigskip
\noindent
Note that we can extend this definition naturally to disconnected graphs by requiring each function in $\text{Lip}(G; h)$ to take a fixed vertex in each connected component of $G$ to $0$. 

\bigskip
\noindent
{\bf{Remark 1.1.2}}. For any $n$-vertex connected graph $G$, $\left\lvert\text{Lip}(G; h)\right\rvert$ is a polynomial of degree $n - 1$ in $\mathbb{Q}[h]$. 

\bigskip
\noindent
{\it{Proof}}. Enumerate the vertices of $G$ as $\{1, 2, \dots, n\}$. Let $P$ be the $\left(n - 1\right)$-dimensional polytope in $\mathbb{R}^{n - 1}$ determined by the hyperplanes $|x_i - x_j| \le 1$ for each edge in $G$ with vertex labels $i$ and $j$, and $x_v = 0$ for a fixed vertex with label $v$. Letting $\mathcal{L} = \mathbb{Z}^{n - 1}$, it's easy to see that $\left\lvert\text{Lip}_{v}(G; h)\right\rvert = \# (hP\ \cap\ \mathcal{L})$ and so $\left\lvert\text{Lip}(G; h)\right\rvert$ is an Ehrhart polynomial [6] of degree $n - 1$. $\blacksquare$

\bigskip
\noindent
Thus when $h$ is large, the growth rate of $\left\lvert\text{Lip}(G; h)\right\rvert$ is governed by the leading coefficient of the associated polynomial. This observation motivates the following definition:

\bigskip
\noindent
{\bf{Definition 1.1.3}}. For any graph $G$ with $n$ vertices and $k < n$ connected components, let $c(G) = \lim_{h \to \infty}\frac{1}{h} \cdot \left\lvert\text{Lip}(G; h)\right\rvert^{1/(n - k)}$. 

\bigskip
\noindent
It's clear due to Remark 1.1.2 that $c(G)$ exists for all graphs $G$. The focus of this paper is to determine exactly or to bound $c(G)$ for a variety of graphs. 

\subsection{Preliminaries}

We begin with basic information about $c(G)$ that will be useful later.

\bigskip
\noindent
{\bf{Theorem 1.2.1}}. For all non-empty graphs $G$, we have that $1 \le c(G) \le 2$. 

\bigskip
\noindent
{\it{Proof}}. Fix $h$ and suppose $v_1, v_2, \dots, v_k \in V(G)$ are our chosen fixed vertices in each of the $k$ connected components of $G$ that $h$-Lipschitz functions will send to $0$. For the lower bound, consider the functions $f: V(G) \rightarrow \mathbb{Z}$ satisfying $f(v_1) = f(v_2) = \dots = f(v_k) = 0$ and $f(w) \in \{0, 1, \dots, h\}$ for all $w \in V(G)\backslash\{v_1, v_2, \dots, v_k\}$. It is clear that all such functions belong to $\text{Lip}(G; h)$ which implies that $\left\lvert\text{Lip}(G; h)\right\rvert \ge (h + 1)^{n - k}$ and thus $c(G) \ge 1$. 

\bigskip
\noindent
For the upper bound, take a maximal spanning forest of $G$ with trees rooted at $v_1, v_2, \dots, v_k$ and perform a depth-first search of each tree. For each vertex $w$, suppose the value of a function $f \in \text{Lip}(G; h)$ at its parent $u$ has already been determined. Then we must have that $f(w) \in \{f(u) - h, f(u) - h + 1, \dots, f(u) + h\}$ so there are at most $2h + 1$ possible values that each $f(w)$ can attain. Therefore we have $\left\lvert\text{Lip}(G; h)\right\rvert \le (2h + 1)^{n - k}$ and thus $c(G) \le 2$. $\blacksquare$

\bigskip
\noindent
Intuitively, dense graphs should have $c(G)$ closer to $1$, while sparse graphs should have $c(G)$ closer to $2$. We formalize this intuition using the most extreme examples of dense graphs (complete graphs) and connected sparse graphs (trees) below:

\bigskip
\noindent
{\bf{Theorem 1.2.2}}. $\left\lvert\text{Lip}(T; h)\right\rvert = (2h + 1)^{n - 1}$ for all trees $T$ with $n$ vertices.

\bigskip
\noindent
{\it{Proof}}. This follows immediately from the proof of the upper bound in Theorem 1.2.1, where we have equality. Note that this implies that $c(T) = 2$. $\blacksquare$

\bigskip
\noindent
{\bf{Corollary 1.2.3}}. If we are counting the number of $h$-Lipschitz functions $f$ on a connected $n$-vertex graph $G$ and the values of $f$ on $j$ vertices have already been determined, there are at most $(2h + 1)^{n - j}$ ways to assign the remaining values of $f$. 

\bigskip
\noindent
We will make extensive use of Corollary 1.2.3 in Section 3 of the paper. 

\bigskip
\noindent
{\bf{Theorem 1.2.4}}. $\left\lvert\text{Lip}(K_n; h)\right\rvert = (h + 1)^n - h^n$ for the $n$-vertex complete graph $K_n$.

\bigskip
\noindent
{\it{Proof}}. For each $f \in \text{Lip}(K_n; h)$, consider the function $g = f - \min_{w \in V(K_n)}f(w)$. Then $g$ is a function from $V(K_n)$ to $\mathbb{N}$ that includes $0$ in its range and satisfies $|g(u) - g(w)| \le h$ for each pair of vertices $u, w \in K_n$. It is clear that there is a bijective mapping between functions $f$ and $g$, so it suffices to count the number of such functions $g$. 

\bigskip
\noindent
Since $0$ is in the range of $g$, we must have that $g(w) \in \{0, 1, \dots, h\}$ for all $w \in V(K_n)$. There are then $(h + 1)^n$ ways in which to construct the function $g$. However $h^n$ of those do not include $0$ in their range, hence the total number of such functions is $(h + 1)^n - h^n$ as desired. Note that this implies that $c(K_n) = \sqrt[n - 1]{n}$. $\blacksquare$

\subsection{Main Results}

Below we present our main results, as well as discuss the organization of the paper:

\bigskip
\noindent
Let $L_{m,n}$ denote the $m \times n$ grid graph (the graph Cartesian product of two paths of length $m$ and $n$), and let $G(n, p)$ denote the Erd\H{o}s-Renyi random graph on $n$ vertices with edge probability $p$. Additionally, let $\alpha \approx 1.16234$ be the largest solution to $\tan{(1/x)} = x$. Our main results are as follows:

\bigskip
\noindent
\begin{itemize}
\item
$\lim_{n \to \infty}c(L_{2, n}) = \alpha\sqrt{2} \approx 1.6437$
\item
$1.351 \approx \alpha^2 \le \lim_{n \to \infty}c(L_{n, n}) \le \frac{1}{\arctan(3/4)} \approx 1.554$
\item
$\lim_{n \to \infty}\mathbf{Pr}\left[1 + \frac{1}{2d} + O\left(\frac{1}{d^2}\right) \le c\left(G\left(n, d/n\right)\right) \le 1 + \frac{4\ln^2{d}}{d} + O\left(\frac{1}{d}\right)\right] = 1$ 
\end{itemize}

\bigskip
\noindent
The paper is organized as follows: In Section $2$ we study $h$-Lipschitz functions on grid graphs, and in Section $3$ we study $h$-Lipschitz functions on sparse random graphs. 

\section{Grid Graphs}

Grid graphs are a particularly interesting domain in which to study $h$-Lipschitz functions, in part due to their natural connection with surface models in statistical mechanics. These graphs have unique characteristics which prevent us from using some of the techniques we apply later for random graphs; specifically, the diameter of grid graphs is polynomial in number of vertices, rather than logarithmic.

\subsection{Upper-Bounding $\lim_{n \to \infty}c(L_{n, n})$}

Our main idea for upper bounding $c(L_{n, n})$ is to upper bound the number of ways to construct an $h$-Lipschitz function $f$ on each row of the grid, given that the values of $f$ have already been determined for the row above. We first utilize a beautiful lemma that lets us assume that $f$ attains the value $0$ on each vertex in the row above, which makes it easy to write a recurrence for the number of ways we can construct $f$ on our given row. Taking $h \to \infty$ then allows us to convert this recurrence into a solvable differential equation. 

\bigskip
\noindent
We begin our exposition with this lemma:

\bigskip
\noindent
{\bf{Lemma 2.1.1}}. Given two centrally-symmetric convex sets $A, B$ in $\mathbb{R}^n$ where $A\cap(B + v)$ is compact for any vector $v \in \mathbb{R}^n$, let $f:\mathbb{R}^n \rightarrow \mathbb{R}$ be the function $f(v) = \text{Vol}(A\cap(B + v))$, where $\text{Vol}$ denotes $n$-dimensional volume. Then $f(v)$ is maximized at $v = 0$. 

\bigskip
\noindent
{\it{Proof}}. Fix $v$, and note that by symmetry $f(v) = f(-v)$. Consider any two $x \in A \cap (B + v)$ and $y \in A \cap (B - v)$. Since $x, y \in A$ by convexity we have that $\frac{x + y}{2} \in A$. Similarly, since $x - v, y + v \in B$ we have that $\frac{(x - v) + (y + v)}{2} \in B$, so $\frac{x + y}{2} \in A \cap B$. Therefore $\frac{1}{2}\left(A \cap (B + v)\right) + \frac{1}{2}\left(A \cap (B - v)\right) \subseteq A \cap B$. 

\bigskip
\noindent
Now, by the Brunn-Minkowski Theorem [7] we have that
\begin{align*}
 f(0) &= \text{Vol}(A \cap B) \\
&\ge \text{Vol}\left(\frac{1}{2}\left(A \cap (B + v)\right) + \frac{1}{2}\left(A \cap (B - v)\right)\right) \\
&\ge \left(\text{Vol}\left(\frac{1}{2}\left(A \cap (B + v)\right)\right)^{1/n} + \text{Vol}\left(\frac{1}{2}\left(A \cap (B - v)\right)\right)^{1/n}\right)^n \\
&= \left(\frac{1}{2} \cdot f(v)^{1/n} + \frac{1}{2} \cdot f(-v)^{1/n}\right)^n \\
&=f(v)
\end{align*}
as desired. $\blacksquare$

\bigskip
\noindent
This lemma is quite useful for upper-bounding $c(G)$ because it will let us show that the scenario in which the number of $h$-Lipschitz functions on a given set of vertices is (roughly) maximized is when $f$ attains the value of $0$ on the neighboring vertices. We formalize this below:

\bigskip
\noindent
{\bf{Corollary 2.1.2}}. Let $G$ be an $n$-vertex connected graph and let $T = \{v_1, v_2, \dots, v_m\} \subset V(G)$ be a proper subset of vertices in $G$. For any $w = (w_1, w_2, \dots, w_m) \in \mathbb{Z}^{m}$ with $w_1 = 0$, let $V_w$ denote the set of $h$-Lipschitz functions $f$ (that fix $f(v_1) = 0$) on $G$ satisfying $f(v_i) = w_i$ for all $i$. Then $|V_0| \ge \left(1 - O(h^{-1})\right)\max_{w}|V_w|$

\bigskip
\noindent
{\it{Proof}}. Let $V(G) = \{v_1, v_2, \dots, v_n\}$ (the first $m$ of these vertices belong to $T$). Additionally, let $A$ be the set of points $(x_1, \dots, x_n) \in \mathbb{R}^{n}$ satisfying $|x_i - x_j| \le h$ for all $\{v_i, v_j\} \in E(G)$ and let $B$ be the set of points in $\mathbb{R}^{n}$ satisfying $x_i = 0$ for all $i \in \{1, 2, \dots, m\}$. Finally, let $w^* = (w_1, w_2, \dots, w_m, 0, 0, \dots, 0) \in \mathbb{Z}^n$. It is clear that there is a one-to-one correspondence between lattice points in $A \cap (B + w^*)$ and $h$-Lipschitz functions in $V_w$. 

\bigskip
\noindent
Letting $\text{Vol}$ denote $(n - m)$-dimensional volume, it is geometrically evident [8] that 

$$\# \left(\left(A \cap (B + w^*)\right) \cap \mathbb{Z}^{n}\right) - \text{Vol}(A \cap (B + w^*)) = O(h^{n - m - 1})$$ 

\bigskip
\noindent
Because $A$ and $B$ are both convex and centrally-symmetric, an application of Lemma 2.1.1 then immediately implies that $|V_0| \ge \left(1 - O(h^{-1})\right)|V_w|$ as desired. $\blacksquare$

\bigskip
\noindent
{\bf{Corollary 2.1.3}}. If the subgraph $G[T]$ induced by the vertices in $T$ is connected, then $|V_0| \ge \left(1 - O(h^{-1})\right)(2h + 1)^{1 - |T|}\left\lvert\text{Lip}(G; h)\right\rvert$.

\bigskip
\noindent
{\it{Proof}}. Because $G[T]$ is connected, the argument from the proof of the upper bound in Theorem 1.2.1 implies that the number of possible vectors $w$ satisfying the Lipschitz condition is $(2h + 1)^{|T| - 1}$. The result then follows immediately from Corollary 2.1.2. $\blacksquare$

\bigskip
\noindent
Using the idea from Corollary 2.1.2, we can prove our first main result:

\bigskip
\noindent
{\bf{Theorem 2.1.4}}. $\lim_{n \to \infty}c(L_{n, n}) \le \frac{1}{\arctan(3/4)} \approx 1.554$. 

\bigskip
\noindent
{\it{Proof}}. Fix $n$, and suppose we have some $h$-Lipschitz function $f$ on our grid graph $L_{n,n}$. Consider some row (other than the top row) of the grid, and denote by $v_1, v_2, \dots, v_n$ the vertices (in order) of this row. Furthermore, fix the values (in order) of $f$ on the vertices in the row above this row as  $w = (w_1, w_2, \dots, w_n) \in \mathbb{Z}^n$.

\bigskip
\noindent
Let $\beta = \arctan(3/4)^{-1} \approx 1.554$. The key idea will be to prove that there are $O\left((\beta{h} + o(h))^n\right)$ ways to assign values of $f(v_i)$ that satisfy the Lipschitz condition - this will imply the theorem, since we can build $f$ one row at a time and bound the number of possibilities each time. 

\bigskip
\noindent
By Corollary 2.1.2, for bounding purposes we can assume that $w = 0$. Under this assumption, let $a(r, s)$ denote the number of ways to assign values to each of $f(v_1), f(v_2), \dots, f(v_r)$ such that $f(v_r) = s$. We immediately obtain the recurrence $a(r, s) = \sum_{i = s - h}^{h}a(r - 1, i)$ for $s \ge 0$ (note that by symmetry $a(r, s) = a(r, -s)$). Letting $a(r)$ be the column vector $\left[a(r, -h), a(r, -h + 1), \dots, a(r, h)\right]^T$ this recurrence is equivalent to the equation:

$$ M_h \cdot a(r - 1) = a(r) $$

\bigskip
\noindent
where $M_h$ is a binary $(2h + 1) \times (2h + 1)$ matrix with entries satisfying $m_{ij} = 1$ when $|i - j| \le h$ and $m_{ij} = 0$ otherwise. 

\bigskip
\noindent
Clearly, to determine the growth rate of $a(r, s)$ it suffices to determine the largest eigenvalue of $M_h$. As $h \to \infty$, $M_h$ becomes a more and more granular discrete approximation of the linear operator $M_{\infty}$ on $\mathcal{L}^2([-1, 1])$ that satisfies (for $x \ge 0$)

$$ M_{\infty} \cdot g(x) = \int_{x - 1}^{1}g(t)\ dt $$

\bigskip
\noindent
Now let $b(x)$ be the eigenfunction of $M_{\infty}$ corresponding to its largest eigenvalue - because $b$ is symmetric about $0$, we may restrict our attention to the interval $[0, 1]$. It is clear that $b(x)$ satisfies the following differential equation for some eigenvalue $\lambda$:

\begin{align}
\lambda \cdot b(x) &= \int_0^{1}b(t)\ dt + \int_0^{1 - x}b(t)\ dt
\end{align}

\noindent
Differentiating both sides we obtain 

\begin{align}
\lambda \cdot b'(x) &= -b(1 - x)
\end{align}

\noindent
and differentiating again we obtain the second-order linear differential equation 

\begin{align}
\lambda^2 \cdot b''(x) &= -b(x)
\end{align}

\noindent
The solution to equation (3) is then given by $b(x) = c_1\cos\left(\frac{x}{\lambda}\right) + c_2\sin\left(\frac{x}{\lambda}\right)$ for some $c_1, c_2 \in \mathbb{R}$. Plugging this solution in to the equation (2) and using the cosine and sine sum formulas, we match terms and obtain

$$ c_1\cos\left(\frac{1}{\lambda}\right) + c_2\sin\left(\frac{1}{\lambda}\right) = c_2 $$

\bigskip
\noindent
Additionally, by letting $x = 0$ and $x = 1$ in equation (1) we see that $b(1) = 2b(0)$, or 

$$ c_1\cos\left(\frac{1}{\lambda}\right) + c_2\sin\left(\frac{1}{\lambda}\right) = 2c_1 $$

\bigskip
\noindent
Thus $\frac{1}{\lambda}$ is the solution to $\cos(x) + 2\sin(x) = 2$ which we can easily solve to find $\lambda = \frac{1}{\arctan(3/4)} = \beta$. 

\bigskip
\noindent
Returning to discrete space, this implies that the largest eigenvalue of $M_h$ is given by $\beta{h} + o(h)$, so $a(r, s) = O\left((\beta{h} + o(h))^{r - 1}\right)$ for all $r$ and $s$. Taking $r = n$ we have

$$ \sum_{s = -h}^{h}a(n, s) = O\left((\beta{h} + o(h))^n\right) $$

\bigskip
\noindent
as desired. The theorem then immediately follows from an application of Corollary 2.1.2. $\blacksquare$

\subsection{Lower-Bounding $\lim_{n \to \infty}c(L_{n, n})$}

We first take a brief detour, and devote our attention to $L_{2,n}$. We will find $\lim_{n \to \infty}c(L_{2, n})$ using the same techniques as in the proof of Theorem 2.1.4, although the resulting differential equation will be slightly more complex. 

\bigskip
\noindent
Our method of proof for the $L_{n,n}$ case will then be to analyze $\text{Lip}(L_{2,n}; h)$ and use our findings to lower bound the average number of ways we can assign values to an $h$-Lipschitz function $f$ on a row of the grid $L_{n,n}$ given that the values on the row below have already been determined. 

\bigskip
\noindent
{\bf{Theorem 2.2.1}}. $\lim_{n \to \infty}c(L_{2, n}) =\alpha\sqrt{2} \approx 1.6437$, where  $\alpha$ is the largest solution to $x = \tan\left(\frac{1}{x}\right)$.

\bigskip
\noindent
{\it{Proof}}. Fix $n$, and suppose we have some $h$-Lipschitz function $f$ on our grid graph $L_{2,n}$ that takes the vertex in the top left corner of the grid to $0$. Suppose we assign values to $f$ one column at a time, and let $a(r, s)$ denote the number of ways to assign values to the first $r$ columns of $L_{2,n}$ such that if $v_{r,1}, v_{r,2}$ are the top and bottom vertices respectively in the $r$th column of $L_{2,n}$, then $f(v_{r,1}) - f(v_{r,2}) = s$. Note that $a(1, s) = 1$ for all $s \in \{-h, -h + 1, \dots, h\}$. 

\bigskip
\noindent
We immediately obtain the recurrence 

$$ a(r, s) = \sum_{i = -h}^{h}a(r - 1, i) \cdot (2h + 1 - |s - i|) $$  

\bigskip
\noindent
Letting $a(r)$ be the column vector $h^{-r + 1} \cdot \left[a(r, -h), a(r, -h + 1), \dots, a(r, h)\right]^T$ this recurrence is equivalent to the equation:

$$ M_h \cdot a(r - 1) = a(r) $$

\bigskip
\noindent
where $M_h$ is a $(2h + 1) \times (2h + 1)$ matrix with entries satisfying $m_{ij} = 2 - \frac{|i - j| - 1}{h}$. 

\bigskip
\noindent
Clearly, to determine the growth rate of $h^{-r + 1} \cdot a(r, s)$ it suffices to determine the largest eigenvalue of $M_h$. By taking $h \to \infty$, we can (as in the proof of Theorem 2.1.4) consider the corresponding linear operator $M_{\infty}$ on $\mathcal{L}^2([-1, 1])$. Let $b(x)$ be the eigenfunction of $M_{\infty}$ corresponding to its largest eigenvalue $\lambda$ - it is clear that $b(x)$ satisfies:

\begin{align}
\lambda \cdot b(x) &= \int_{-1}^{x}(2 + t - x)b(t)\ dt + \int_{x}^{1}(2 - t + x)b(t)\ dt
\end{align}

\noindent
Note that $b(x) = b(-x)$ so we will restrict our attention to $x \ge 0$. Differentiating both sides with respect to $x$ we obtain 

\begin{align}
\lambda \cdot b'(x) &= \left(2b'(x) -\int_{-1}^{x}b(t)\ dt\right) + \left(-2b'(x) + \int_{x}^{1}b(t)\ dt\right) = -2\int_0^{x}b(t)\ dt
\end{align}

\noindent
and differentiating again we obtain the second-order linear differential equation 

\begin{align}
\lambda \cdot b''(x) &= -2b(x)
\end{align}

\noindent
Let $\gamma = \sqrt{2/\lambda}$. Once again the solution to equation (6) is then given by $b(x) = c_1\cos\left(\gamma{x}\right) + c_2\sin\left(\gamma{x}\right)$ for some $c_1, c_2 \in \mathbb{R}$. But note that from equation (5) we have that $b'(0) = 0$, so $c_2 = 0$ and $b(x) = c_1\cos\left(\gamma{x}\right)$, and we may assume without loss of generality that $c_1 = 1$. Finally, plugging this expression in to the equation (4) we obtain 

$$ \lambda \cdot b(x) = \frac{2}{\gamma^2} \cdot b(x) + \frac{2}{\gamma} \cdot \sin(\gamma) - \frac{2}{\gamma^2}\cos(\gamma) $$

\bigskip
\noindent
But since $\lambda = \frac{2}{\gamma^2}$ after canceling we have that $\gamma$ is the smallest positive solution to $\tan(\gamma) = \frac{1}{\gamma}$. Thus $\gamma = \frac{1}{\alpha}$ and $\lambda = 2\alpha^2$. 

\bigskip
\noindent
Returning to discrete space, this implies that the largest eigenvalue of $M_h$ is given by $2\alpha^2{h} + o(h)$, so $h^{-r + 1} \cdot a(r, s) = \Omega\left(\left(2\alpha^2\cdot{h} + o(h)\right)^{r - 1}\right)$ for all $s$. Taking $r = n$ we have

$$ \sum_{s = -h}^{h}a(n, s) = \Omega\left((\alpha\sqrt{2}h + o(h))^{2n - 1}\right) $$

\bigskip
\noindent
which immediately implies that  $\lim_{n \to \infty}c(L_{2, n}) = \alpha\sqrt{2}$ as desired. $\blacksquare$

\bigskip
\noindent
We can use this result about $L_{2,n}$ to obtain an upper bound on  $\lim_{n \to \infty}c(L_{n, n})$. In particular, we can show:

\bigskip
\noindent
{\bf{Theorem 2.2.2}}. $\lim_{n \to \infty}c(L_{n, n}) \ge \alpha^2 \approx  1.35103$.

\bigskip
\noindent
{\it{Proof}}. Fix $n$ and let $S \subset \mathbb{Z}^n$ be the set of all vectors $(0, x_2, \dots, x_n) \in \mathbb{Z}^n$ such that $|x_2| \le h$ and $|x_{i + 1} - x_{i}| \le h$ for all $i \ge 2$. For simplicity, we will focus on constructing $h$-Lipschitz functions $f$ that attain $0$ on each vertex of the left-most column of the $L_{n,n}$ grid. 

\bigskip
\noindent
Call two vectors $x, y \in S$ mutually {\it{close}} if $\|x - y\|_{\infty} \le h$ and, let $S(x) \subset S$ denote the set of vectors $y \in S$ such that $y$ is close to $x$. We will build an $h$-Lipschitz function $f$ on the vertices of the $L_{n,n}$ one row at a time, and count how many such functions exist.

\bigskip
\noindent
Let $a(k, x)$ denote the number of $h$-Lipschitz functions $f$ on $L_{k,n}$ where $x \in S$ corresponds to the values of $f$ of the top row, and where $f$ attains $0$ on each vertex in the left-most column of $L_{k,n}$. Then it's clear that $a(k, x) = \sum_{y \in S(x)}a(k - 1, y)$. 

\bigskip
\noindent
Now, enumerate the $(2h + 1)^{n - 1}$ elements of $S$. In the spirit of the proofs of Theorems 2.1.4 and 2.2.1, let $M_h$ be the $(2h + 1)^{n - 1} \times (2h + 1)^{n - 1}$ binary matrix with entries satisfying $m_{ij} = 1$ if elements $i$ and $j$ of $S$ are close and $0$ otherwise. Based on the recurrence in the preceding paragraph, the growth rate of $a(k, x)$ corresponds to the largest eigenvalue of $M_h$. 

\bigskip
\noindent
In the proof of Theorem 2.2.1, we established that $\left|\text{Lip}_v(L_{2,n})\right| = \Omega\left((\alpha\sqrt{2} \cdot h + o(h))^{2n - 1}\right)$ where $v$ denotes the vertex in the bottom left-hand corner of the $L_{2,n}$ grid. Corollary 2.1.3 then implies that the number of $h$-Lipschitz functions on $L_{2,n}$ that attain $0$ on {\it{both}} vertices in the left-most column of $L_{2,n}$ is least $\left(1 - O(h^{-1})\right)(2h + 1)^{-1}\left|\text{Lip}_v(L_{2,n})\right| = \omega\left((\alpha\sqrt{2} \cdot h + o(h))^{2n - 2}\right)$. 

\bigskip
\noindent
But note also that the number of such functions $f$ is given by $\sum_{x \in S}|S(x)|$, which can be seen by picking any vector in $x \in S$ to correspond to the values of $f$ on the bottom row of the grid, and then any vector in $S(x)$ to correspond to the values of $f$ on the row above. 

\bigskip
\noindent
Therefore we have that ${\bf{1}}^T \cdot M_h \cdot {\bf{1}} = \sum_{x \in S}|S(x)| = \omega\left((\alpha\sqrt{2} \cdot h + o(h))^{2n - 2}\right)$. But then we know that the largest eigenvalue $\lambda$ of $M_h$ satisfies 

$$ \lambda \ge \frac{{\bf{1}}^T \cdot M_h \cdot {\bf{1}}}{{\bf{1}}^T \cdot {\bf{1}}} = \omega\left(\left(\alpha^2{h} + o(h)\right)^{n - 1}\right)$$

\bigskip
\noindent
which immediately implies the result. $\blacksquare$

\subsection{Improving Bounds on $\lim_{n \to \infty}c(L_{n, n})$}

{\bf{Remark 2.3.1}}. We can use the same techniques as in the proofs of Theorems 2.1.4 and 2.2.2 to improve both our upper and lower bound. 

\bigskip
\noindent
Indeed, suppose we knew that the number of $h$-Lipschitz functions on $L_{3,n}$ taking a value of $0$ at each vertex in the top row of the $L_{3,n}$ was given by $\Omega\left((\zeta{h} + o(h))^{2n}\right)$ for some $\zeta \in [1, 2]$. Then repeating the argument from the proof of Theorem 2.1.4, by building an $h$-Lipschitz function on $L_{n,n}$ {\it{two}} rows at a time, we would have that $\lim_{n \to \infty}c(L_{n, n}) \le \zeta$. 

\bigskip
\noindent
Additionally, suppose we knew that $\left|\text{Lip}(L_{3,n}; h)\right| = \Omega\left((\psi{h} + o(h))^{3n - 1}\right)$ for some $\psi \in [1, 2]$. Then, re-using $M_h$ and $\lambda$ and repeating the argument from the proof of Theorem 2.2.2, we have that ${\bf{1}}^T \cdot M_h^2 \cdot {\bf{1}}$ is at least $\left(1 - O(h^{-1})\right)(2h + 1)^{-2}\left|\text{Lip}(L_{3,n}; h)\right|$. Therefore

$$ \lambda^2 \ge \frac{{\bf{1}}^T \cdot M_h^2 \cdot {\bf{1}}}{{\bf{1}}^T \cdot {\bf{1}}} = \omega\left(\left(\frac{\psi^3{h^2}}{2} + o(h^2)\right)^{n - 1}\right) $$

\bigskip
\noindent
which implies a lower bound of $\lim_{n \to \infty}c(L_{n, n}) \ge \frac{\psi^{3/2}}{\sqrt{2}}$

\bigskip
\noindent
Calculating $\zeta$ is quite similar to calculating $\lim_{n \to \infty}c(L_{2, n})$, although it requires solving the following more complicated differential equation:

$$ \zeta^2 \cdot b(x, y) = \iint_{R}b(s, t)\ ds\ dt $$

\bigskip
\noindent
where $R$ is the region in $\mathbb{R}^2$ determined by $|s|, |t| \le 1$ and $|x - s| \le 1$ and $|x + y - s - t| \le 1$.

\bigskip
\noindent
Calculating $\psi$ is more difficult still, with corresponding differential equation:

$$ \psi^3 \cdot b(x, y) = \int_{-1}^{1}\iint_{R_k}b(s, t)\ ds\ dt\ dk $$

\bigskip
\noindent
where $R_k$ is the region in $\mathbb{R}^2$ determined by $|s|, |t| \le 1$ and $|x - s - k| \le 1$ and $|x + y - s - t - k| \le 1$. 

\bigskip
\noindent
We leave precise calculation of $\zeta$ and $\psi$ for future work, although a computer-assisted calculation suggests that $\zeta \approx 1.4895$ and $\psi \approx 1.553$. Applying these approximations results in lower and upper bounds on $\lim_{n \to \infty}c(L_{n, n})$ of $1.3685$ and $1.4895$ respectively.

\section{Random Graphs}

Random graphs are another natural class of graphs on which to study $c(G)$. Because these graphs have much less defined structure than grid graphs, our methods in this section are understandably more probabilistic in nature. 

\subsection{Lower-Bounding $\lim_{n \to \infty}c\left(G\left(n, d/n\right)\right)$}

We begin with our lower bound, which utilizes the Los\'avz Local Lemma (LLL) [9]. The main idea is to carefully construct a large number of potential functions $f$, and prove that a sufficiently large proportion of them satisfy the $h$-Lipschitz conditions. 

\bigskip
\noindent
{\bf{Theorem 3.1.1}}. $\lim_{n \to \infty}\mathbf{Pr}\left[c\left(G\left(n, d/n\right)\right) \ge 1 + \frac{1}{2d} + O\left(\frac{1}{d^2}\right)\right] = 1$.

\bigskip
\noindent
{\it{Proof}}. Let $G = G\left(n, d/n\right)$ for some fixed $n$. It is well-known [10] that as $n$ increases, the degree distribution of $G$ approaches a $\text{Poisson}(d)$ distribution. Using the standard Poisson tail bound

$$ \mathbf{Pr}(X \ge d + x) \le e^{-\frac{x^2}{2(x + d)}} $$

\bigskip
\noindent
we have that the proportion of vertices in $G$ of degree at least $2d$ is with high probability at most $2e^{-d/4}$. 

\bigskip
\noindent
Now, consider a function $f$ on the vertices of $G$. For each vertex $v \in V(G)$ with degree less than $2d$, assign to $f(v)$ uniformly at random an integer in $\{0, 1, \dots, \lfloor{(1 + c)h}\rfloor\}$ for some constant $c > 0$ to be determined later. For all other vertices $v$, assign to $f(v)$ uniformly at random an integer in $\{\lceil{ch}\rceil, \lceil{ch}\rceil + 1, \dots, h\}$. Using the bound proven above as well as the fact that $(1 + c)(1 - c) < 1$, with high probability there are at least $((1 + c)(1 - c)^{4e^{-d/4}}h + o(h))^{n}$ such functions. 

\bigskip
\noindent
Now, we will lower-bound the probability that such a function $f$ is $h$-Lipschitz. First, we will modify $f$ slightly so that it satisfies the condition that a fixed vertex in each connected component in $G$ is sent to $0$. Suppose there are $k$ connected components in $G$; after our initial construction of $f$, we level shift $f$ in each connected component so that this condition is satisfied. Clearly this transformation is at most $(1 + c)^kh^k$-to-one (which can be seen by considering the original value of $f$ on each fixed vertex in each connected component of $G$). 

\bigskip
\noindent
Furthermore, it is well-known [10] that with high probability, the giant component of $G$ contains at least $\left(1 - \frac{x}{d} + o(1)\right)n$ vertices, where $x  = -W\left(-de^{-d}\right) \approx de^{-d}$ and $W$ is the Lambert $W$-function. Thus with high probability we have that $k \le e^{-d/4}n$ for sufficiently large $d$, so that the number of distinct transformed functions $f$ is at least $\left((1 + c)(1 - c)^{5e^{-d/4}} + o(1)\right)^n \cdot h^{n - k}$, where we once again use the fact that $(1 + c)(1 - c) < 1$. With the $0$ condition accounted for, it remains to lower bound the probability that $|f(a) - f(b)| \le h$ for all $\{a, b\} \in E(G)$. 

\bigskip
\noindent
Enumerate each edge in $G$, and let $A_i$ be the event that edge $i$ fails the Lipschitz condition. With high probability, there are $\frac{dn}{2} + o(n)$ such events. It is easy to see that any edge incident to a vertex of degree at least $2d$ cannot fail the Lipschitz condition, and for all other edges $i$ we have that $\mathbf{Pr}(A_i) = \frac{c^2}{(1 + c)^2} + o(1)$. Furthermore, each of these $A_i$ is only dependent on events $A_j$ where edge $j$ is incident to edge $i$, of which there are at most $4d - 4$. Thus, letting $h \to \infty$ so that we can ignore the $o(1)$ terms, by LLL we have that for any $x > 0$

$$ \frac{c^2}{(1 + c)^2} \le x(1 - x)^{4d - 4} \Longrightarrow \mathbf{Pr}\left(\bigcup_{i}\bar{A_i}\right) \ge (1 - x)^{\frac{dn}{2} + o(n)} $$

\bigskip
\noindent
Let $x = \frac{1}{d(d - 1)} - \frac{1}{2d^2(d - 1)^2}$ and $c = \frac{1}{d}\sqrt{1 - \frac{4}{d}}$. Noting that $1 - r + \frac{r^2}{2} \ge e^{-r} \ge 1 - r$ for all $r \ge 0$ we have that $(1 - x)^{4d - 4} \ge e^{-4/d} \ge 1 - \frac{4}{d}$. Furthermore, it is easy to compute that $x \ge \frac{1}{d^2}$. Therefore

$$ \frac{c^2}{(1 + c)^2} \le c^2 = \frac{1}{d^2}\left(1 - \frac{4}{d}\right) \le x(1 - x)^{4d - 4} $$

\bigskip
\noindent
Thus for this choice of $x$ and $c$ we have 

$$ \mathbf{Pr}\left(\bigcup_{i}\bar{A_i}\right) \ge (1 - x)^{\frac{dn}{2} + o(n)} \ge e^{-\frac{1}{d(d - 1)}\left(\frac{dn}{2} + o(n)\right)} \ge \left(1 - \frac{1}{d - 1}\right)^{n/2 + o(n)}$$

\bigskip
\noindent
Multiplying this lower bound on the probability that a randomly selected $f$ is $h$-Lipschitz by the number of potential functions $f$ and noting that

$$ \left(1 +  \frac{1}{d}\sqrt{1 - \frac{4}{d}}\right)\left(1 -  \frac{1}{d}\sqrt{1 - \frac{4}{d}}\right)^{5e^{-d/4}}\sqrt{1 - \frac{1}{d - 1}} = 1 + \frac{1}{2d} + O\left(\frac{1}{d^2}\right) $$

\bigskip
\noindent
then completes the proof. $\blacksquare$

\subsection{Upper-Bounding $\lim_{n \to \infty}c\left(G\left(n, d/n\right)\right)$}

We continue with the upper bound. The main idea follows an argument made in [1], which is that for any $h$-Lipschitz function $f$ on a graph $G$, the set of vertices $\{v \in V(G)\ |\ f(v) \le j\}$ and the set of vertices $\{v \in V(G)\ |\ f(v) > j + h\}$ can have no edges between them. We will show that in a random graph, with high probability there do not exist two large sets of vertices with no edges between them. This will prove that any $h$-Lipschitz function on a random graph must be ``flat" (i.e. take only $h + 1$ distinct values) on a large portion of the graph. 

\bigskip
\noindent
We will then provide a sketch of an argument showing that the function $f$ is extremely constrained on the majority of vertices outside the ``flat" portion of the graph, which improves the bound further for most random graphs.

\bigskip
\noindent
We begin with a helpful lemma:

\bigskip
\noindent
{\bf{Lemma 3.2.1}}. Let $G = G(n, d/n)$ and let $\alpha = \frac{2\ln{d}}{d}$. Then as $n \to \infty$, with high probability there do not exist disjoint vertex sets $A, B \subset V(G)$ satisfying $\min\left(|A|, |B|\right) \ge \alpha{n}$ with no edge between them. 

\bigskip
\noindent
{\it{Proof}}. let $H$ denote the standard entropy function $H(p) = -p\log_2(p) - (1 - p)\log_2(1 - p)$ and let $X$ be a random variable representing the number of pairs of disjoint sets $A, B \subset V(G)$ satisfying $|A| = |B| = \left\lceil \alpha{n} \right\rceil$ with no edge between them. There at most $\binom{n}{2\left\lceil{\alpha{n}}\right\rceil} \cdot \binom{2\left\lceil{\alpha{n}}\right\rceil}{\left\lceil{\alpha{n}}\right\rceil} < e^{nH(2\alpha)\ln{2}} \cdot e^{2n\alpha\ln{2}}$ pairs of sets $(A, B)$ satisfying the size conditions, and the probability that any such pair contains no edges between them is $\left(1 - \frac{d}{n}\right)^{{\left\lceil{\alpha{n}}\right\rceil}^2} \le e^{-dn\alpha^2}$. Mutliplying these values to upper bound $\mathbb{E}[X]$ and noting that $\mathbf{Pr}(X = 0) \ge 1 - \mathbb{E}[X]$, we see that it suffices to show that $d\alpha^2 - 2\alpha\ln{2} - H(2\alpha)\ln{2} > 0$.

\bigskip
\noindent
Using the inequality $-(1 - x)\ln{(1 - x)} < x$ we have $H(2\alpha) < -2\alpha\log_2\left(\frac{2\alpha}{e}\right)$. Then it's easy to see that

$$ 2\alpha\ln{2} + H(2\alpha)\ln{2} < -2\alpha\ln\left(\frac{\alpha}{e}\right) = \frac{4\ln{d}}{d} \cdot \ln{\left(\frac{ed}{2\ln{d}}\right)} < \frac{4\ln^2{d}}{d} = d\alpha^2 $$

\bigskip
\noindent
which completes the proof. $\blacksquare$

\bigskip
\noindent
Note that a similar lemma was proven in Chapter 6 of [10]. With this lemma in hand, we may now begin the proof of our main result in this subsection:

\bigskip
\noindent
{\bf{Theorem 3.2.2}}. $\lim_{n \to \infty}\mathbf{Pr}\left[c\left(G\left(n, d/n\right)\right) \le 1 + \frac{4\ln^2{d}}{d} + O\left(\frac{1}{d}\right)\right] = 1$.

\bigskip
\noindent
{\it{Proof}}. As before, let $\alpha = \frac{2\ln{d}}{d}$ and let $G = G(n, d/n)$. Furthermore, let $G'$ denote the largest connected component (the giant component) of $G$. Consider an $h$-Lipschitz function $f$ on $G'$. We will first show that there exists $j \in \mathbb{Z}$ such that for a large proportion of vertices $v \in V(G')$, $f(v) \in \{j, j + 1, \dots, j + h\}$. 

\bigskip
\noindent
Indeed let $j$ be the smallest integer for which $|\{v \in V(G')\ |\ f(v) \le j\}| \ge \alpha{n}$. Then using Lemma 3.2.1 we have that with high probability $|\{v \in V(G')\ |\ f(v) > j + h\}| < \alpha{n}$, since by the Lipschitz condition there are no edges between these two subsets of $V(G')$. Thus

$$  |\{v \in V(G')\ |\ f(v) \in \{j, j + 1, \dots, j + h\}| > m - 2\alpha{n} $$

\bigskip
\noindent
where $|G'| = m$. Now, let 
\begin{align*}
S &= \{v \in V(G)\ |\ f(v) \in \{j, j + 1, \dots, j + h\}\} \\
S_1 &= \{v \in V(G)\ |\ f(v) < j\} \\
S_2 &= \{v \in V(G)\ |\ f(v) > j + h\} 
\end{align*}

\noindent
and suppose these vertex sets are fixed. We have that $|S_1|,|S_2| < \alpha{n}$ and $|S| > m - 2\alpha{n}$. 

\bigskip
\noindent
Let us now count the number of $h$-Lipschitz functions $f$ satisfying these conditions. Without loss of generality we may assume that $f$ sends a vertex in $S$ to $0$. Following the proof of Theorem 1.2.4, this implies that there are $(h + 1)^{|S|} - h^{|S|} < n(h + 1)^{|S| - 1}$ ways to assign values for $f$ on the vertices in $S$. Once these values have been chosen, using Corollary 1.2.3 there are at most $(2h + 1)^{m - |S|}$ ways to assign the remaining values of $f$. Thus we have that

\begin{align*}
\left\lvert\text{Lip}(G'; h)\right\rvert &\le \sum_{\substack{S \subseteq [m] \\ |S| > m - 2\alpha{n}}}n(h + 1)^{|S| - 1}(2h + 1)^{m - |S|} \\
&\le n(h + 1)^{m - 1} \cdot \sum_{i \le 2\alpha{n}}\binom{m}{i} \cdot 2^{2\alpha{n}} \\
&\le n(h + 1)^{m - 1} \cdot 2^{H(2\alpha)n} \cdot 2^{2\alpha{n}} \\
& \le n(h + 1)^{m - 1} \cdot e^{d\alpha^2{n}}
\end{align*}

\bigskip
\noindent
where the last inequality follows from the fact that $H(2\alpha)\ln{2} + 2\alpha\ln{2} < d\alpha^2$ as shown in the proof of Lemma 3.2.1. If $G$ has $k$ connected components, Corollary 1.2.3 implies that 

$$ \left\lvert\text{Lip}(G; h)\right\rvert \le \left\lvert\text{Lip}(G'; h)\right\rvert \cdot (2h + 1)^{n - m - k + 1} $$

\bigskip
\noindent
But as discussed in the proof of Theorem 3.1.1 we have with high probability that $m \ge (1 - e^{-d/4})n$, which immediately implies a bound of 

$$ 2^{e^{-d/4}} \cdot e^{d\alpha^2 \cdot (1 - e^{-d/4})^{-1}} = 1 + d\alpha^2 + O\left(\frac{1}{d}\right) = 1 + \frac{4\ln^2{d}}{d} + O\left(\frac{1}{d}\right) $$
  
\bigskip
\noindent
as desired. $\blacksquare$

\bigskip
\noindent
{\bf{Remark 3.2.3}}. We will now provide a sketch of an argument that shows for {\it{most}} random graphs, we can actually improve our upper bound to $1 + O\left(\frac{\ln{d}}{d}\right)$. In particular, that there exist fixed constants $b, B> 1$ such that for all sufficiently large $d$, $\lim_{n \to \infty}\mathbf{Pr}\left[c\left(G\left(n, d/n\right)\right) \le 1 + \frac{b\ln{d}}{d} + O\left(\frac{1}{d}\right)\right] \ge 1 - B^{-d/\ln{d}}$.

\bigskip
\noindent
{\it{Sketch.}} As before, let $\alpha = \frac{2\ln{d}}{d}$ and let $G = G(n, d/n)$. For simplicity, we will assume $G$ is connected. Once again we let $S$ denote the set of vertices for which an $h$-Lipschitz function on $G$ is ``flat;" from the proof of Theorem 3.2.2, we have that $|S| \ge (1 - 2\alpha)n$. In our proof of Theorem 3.2.2, for vertices in $V(G)\backslash S$ we bounded the number of possible values of $f$ on each of these vertices by $2h + 1$. It turns out that this bound is incredibly weak; we claim now that for the majority of these vertices, the number of possible values of $f$ is actually $O\left(\alpha{h}\right)$. 

\bigskip
\noindent
Indeed, suppose we remove edges in $G$ at random with probability $1 - \Omega\left(\frac{1}{\ln{d}}\right)$. Then by a simple probabilistic alteration argument, for all but a small number of {\it{lonely}} vertices in $V(G) \backslash S$, we can associate with that non-lonely vertex a set of $\Omega\left(\frac{1}{\alpha}\right)$ vertices in $S$ connected to that vertex but connected to no other vertex in $V(G) \backslash S$ (post-edge removal). In particular, use of the Poisson tail bound shows that the number of lonely vertices roughly follows a binomial distribution with mean and variance $\text{exp}\left(-\Omega(\alpha)\right) \cdot n$. 

\bigskip
\noindent
Now, recall that $f$ must take a value between $j$ and $j + h$ for each vertex in $S$. For a given non-lonely vertex in $V(G) \backslash S$, because it has $\Omega\left(\frac{1}{\alpha}\right)$ neighbors in $S$, over all possible values of $f$ on the vertices in $S$ the expected number of possible values of $f$ on this vertex that satisfy the Lipschitz condition is then $\Omega\left(\alpha{h}\right)$. This proves the claim, and by replacing the $2h + 1$ in the proof of Theorem 3.2.2 we find that the claim implies the $1 + O\left(\alpha\right)$ bound, assuming that there are not too many lonely vertices. 

\bigskip
\noindent
The last remaining step is to bound the probability that there is some ``adversarial" $S \subset V(G)$ that (post edge-removal) produces many lonely vertices. Recall that the number of lonely vertices each set $S$ produces roughly follows a binomial distribution with mean and variance $\text{exp}\left(-\Omega(\alpha)\right) \cdot n$ and the number of possible subsets $S$ is at most $2^n$. Furthermore, it is well-known that the expected maximum of $N$ sub-Gaussian random variables each with mean $\mu$ and variance $\sigma^2$ can be bounded by $\mu + \sigma\sqrt{2\log{N}}$. Applying this bound on random variables representing the number of lonely vertices produced by a set $S$ with $N \le 2^n$ and using Markov's inequality then immediately implies the desired result. $\blacksquare$

\bigskip
\noindent
{\bf{Remark 3.2.4.}}  Note that the upper bound in Theorem 3.2.2 is only non-trivial for large $d$. Thus we must utilize a different argument to show that $\lim_{n \to \infty}\mathbf{Pr}\left[c\left(G\left(n, d/n\right)\right) < 2\right] = 1$ for all $d > 1$. 

\bigskip
\noindent
The idea will be to with high probability find a ``ladder-like" structure in $G(n, d/n)$, consisting of two vertex-disjoint paths $P_1$ and $P_2$ each of size $\Omega(n)$, where $\Omega(n)$ distinct vertices in $P_1$ are adjacent to a vertex in $P_2$. We will then repeat arguments similar to those in the proof of Theorem 2.1.4 (two adjacent rows in an $L_{n,n}$ form a similar ladder-like structure) to bound the number of ways we can assign values to an $h$-Lipschitz function $f$ for vertices on path $P_1$. 

\bigskip
\noindent
{\bf{Theorem 3.2.5}}. $\lim_{n \to \infty}\mathbf{Pr}\left[c\left(G\left(n, d/n\right)\right) \le 2 - 2^{-18}\epsilon^5\right] = 1$ for $d = 1 + \epsilon$, where $\epsilon$ is sufficiently small. 

\bigskip
\noindent
{\it{Proof}}. Let $G$ be the graph with $n$ (for sufficiently large $n$) vertices and no edges, and suppose we randomly construct edges in $G$ in two stages; in Stage $1$, for each pair of vertices we construct the edge between them with probability $(1 + \epsilon)/n$, and in Stage $2$ we again attempt to construct those edges with probability $\epsilon/n$ (the same edge may be constructed twice). It is clear that the resulting graph is equivalent to a $G\left(n, (1 + 2\epsilon)/n - \epsilon(1 + \epsilon)/n^2\right) \subseteq G\left(n, (1 + 2\epsilon)/n\right)$ random graph.

\bigskip
\noindent
Let $E_1(G)$ and $E_2(G)$ denote the edges constructed in Stages $1$ and $2$ respectively. By Theorem 3.2 in [11], with high probability there is a path $P$ of length\footnote{All path lengths and vertex counts are integral; we eschew the use of floor functions for the sake of readability} $\epsilon^2n/5$ in $G$ consisting solely of edges in $E_1(G)$. Split this path into two connected halves $P_1$ and $P_2$, with $|P_1|, |P_2| = \epsilon^2n/10$. For each vertex in $P_1$ there are $\epsilon^2n/10$ potential edges between that vertex and a vertex in $P_2$, and each of these edges has probability $\epsilon/n$ of being present in $E_2(G)$. 

\bigskip
\noindent
Let $X$ be a random variable representing the number of vertices in $P_1$ for which there is an edge in $E_2(G)$ connecting that vertex to $P_2$. We have that
 
$$ X \sim \text{Binomial}\left(\frac{\epsilon^2n}{10},  1 - \left(1 - \frac{\epsilon}{n}\right)^{\epsilon^2n/10}\right) \dot\sim\ \text{Binomial}\left(\frac{\epsilon^2n}{10},  \frac{\epsilon^3}{10}\right) $$

\bigskip
\noindent
thus with high probability $X \ge \epsilon^5n/200$. 

\bigskip
\noindent
Let $F$ be a maximal spanning forest of the edge-induced subgraph $G[E_1]$ that contains $P$, and add the edges in $E_2(G)$ (so that $F$ is no longer a forest). We will upper bound $\left\lvert\text{Lip}(F; h)\right\rvert$. 

\bigskip
\noindent
Construct $h$-Lipschitz functions $f$ on $F$ as follows; select a vertex in each of the $k$ connected components of $F$ that $f$ will send to $0$. Then for each edge $\{x, y\} \in E(F) \cap E_1(G)$ choose $f(x) - f(y)$ uniformly at random from $\{-h, -h + 1, \dots, h\}$. This process constructs $(2h + 1)^{n - k}$ functions $f$ with uniform probability, including all valid $h$-Lipschitz functions. We will upper bound the probability that for a given function $f$ constructed in this manner, the edges in $E_2(G)$ satisfy the $h$-Lipschitz condition.

\bigskip
\noindent
Denote by $v_1, v_2, \dots, v_m$ (in order) the vertices in $P_1$ connected to vertices in $P_2$, where with high probability $m \ge \epsilon^5n/200$. Using Corollary 2.1.2, for sufficiently large $h$ we can assume that $f(w) = 0$ for each vertex $w \in P_2$. Then we must have that $|f(v_i)| \le h$ for all $i$. Now, consider every third vertex $v_1, v_4, v_7, \dots$ (there are at least $ \epsilon^5n/600$ such vertices). Then there are at least two vertices between $v_1$ and $v_4$ in path $P_1$, so let these vertices be $u_1, u_2, \dots, u_j$ for some $j \ge 2$. We have that 

$$ |(f(u_1) - f(v_1)) + (f(u_2) - f(u_1)) + \dots + (f(v_3) - f(u_j))| = |f(v_3) - f(v_1)| \le 2h $$  

\bigskip
\noindent
But each of the terms in the left-most summation is uniformly distributed from $\{-h, -h + 1, \dots, h\}$. Clearly the maximum probability that the above inequality is satisfied occurs when $j = 2$, and we can easily compute this probability as approaching $23/24$ as $h \to \infty$. Indeed, geometrically the possibility space is a $2h \times 2h \times 2h$ cube and the space corresponding to failure corresponds to two tri-rectangular tetrahedrons, each with volume $h^3/6$. 

\bigskip
\noindent
We can perform the same calculation with vertices $v_4$ and $v_7$, then $v_7$ and $v_{10}$, and so on. Thus, as $h \to \infty$ the probability that our randomly chosen function $f$ satisfies the Lipschitz condition is at most $(23/24 + o(1))^{\epsilon^5n/600} < \left(1 - 2^{-14}\epsilon^5\right)^{n - k}$. Thus

$$ \left\lvert\text{Lip}(G; h)\right\rvert < \left\lvert\text{Lip}(F; h)\right\rvert \le \left(1 - 2^{-14}\epsilon^5\right)^{n - k} \cdot (2h + 1)^{n - k} $$

\bigskip
\noindent
Replacing $\epsilon$ with $\epsilon/2$ now implies the desired result. $\blacksquare$ 

\section{Acknowledgements}
The authors would like to thank Mark Sellke for pointing us towards Lemma 2.1.1, and Yonah Borns-Weil for his input on the spectral theory of operators. Additionally, we would like to thank Franklyn Wang and Ryan Alweiss for helpful discussions pertaining to our investigations in this paper.

\end{document}